\documentclass[reqno]{amsart}
\usepackage{amscd}
\begin{document}
\def\N{{\mathbb N}}
\def\Z{{\mathbb Z}}
\def\R{{\rm I\!R}}
\def\phi{\varphi}

\def\C{{\mathchoice {\setbox0=\hbox{$\displaystyle\rm C$}\hbox{\hbox to0pt{\kern0.4\wd0\vrule height0.9\ht0\hss}\box0}}{\setbox0=\hbox{$\textstyle\rm C$}\hbox{\hbox to0pt{\kern0.4\wd0\vrule height0.9\ht0\hss}\box0}} {\setbox0=\hbox{$\scriptstyle\rm C$}\hbox{\hbox to0pt{\kern0.4\wd0\vrule height0.9\ht0\hss}\box0}} {\setbox0=\hbox{$\scriptscriptstyle\rm C$}\hbox{\hbox to0pt{\kern0.4\wd0\vrule height0.9\ht0\hss}\box0}}}}

\def\codim{{\rm codim}}
\def\M{{\mathcal M}}
\def\ov{\overline}
\def\m{\mapsto}
\def\rk{{\rm rank}}
\def\id{{\sf id}}
\def\Aut{{\sf Aut}}
\def\CR{{\rm CR}}
\def\crd{\dim_{{\rm CR}}}
\def\crc{{\rm codim_{CR}}}
\emergencystretch9pt
\frenchspacing

\newtheorem{Th}{Theorem}[section]
\newtheorem{Def}{Definition}[section]
\newtheorem{Cor}{Corollary}[section]
\newtheorem{Satz}{Satz}[section]
\newtheorem{Prop}{Proposition}[section]
\newtheorem{Lemma}{Lemma}[section]
\newtheorem{Rem}{Remark}[section]
\newtheorem{Example}{Example}[section]

\begin{abstract}
The topic of the paper is the study of germs
of local holomorphisms $f$ between $C^n$ and $C^{n'}$
such that $f(M)\subset M'$ and $df(T^cM)=T^cM'$
for $M\subset C^n$ and $M'\subset C^{n'}$ generic 
real-analytic CR submanifolds of arbitrary codimensions.
It is proved that for $M$ minimal and $M'$ finitely nondegenerate,
such germs depend analytically on their jets.
As a corollary, an analytic structure on the set of
all germs of this type is obtained.
\end{abstract}

\title[Germs of local automorphisms]{Germs of local automorphisms of real-analytic CR structures and analytic dependence on $k$-jets.}
\author{Dmitri Zaitsev}
\address{Mathematisches Institut, Universit\"at T\"ubingen, 
        72076 T\"ubingen, Germany, 
        E-mail address: dmitri.zaitsev@uni-tuebingen.de}
\thanks{Partially supported by SFB 237 ``Unordnung und gro\ss e Fluktuationen''}
\maketitle

\section{Introduction.}

Let $M\subset\C^n,M'\subset\C^{n'}$ 
be connected locally closed real-analytic submanifolds,
$x\in M, x'\in M'$ be arbitrary points.
The complex tangent subspace $T_xM\cap iT_xM$ 
will be denoted by $T^c_xM$. $M$ is a {\em CR manifold},
if $\dim T^c_xM$ is constant. In this case 
$\crd M:= \dim_{\C}T^cM$ is called the {\em CR dimension} and
$\crc M:= \dim_{\R}TM - \dim_{\R}T^cM$ the {\em CR codimension}.
A CR submanifold $M\subset\C^n$ is called {\em generic},
if $TM+iTM=\C^n$.

Suppose for the moment that $M,M'\subset\C^n$ 
are generic real-analytic {\rm CR} submanifolds
of the same {\rm CR} dimension and the same {\rm CR} codimension.
Baouendi, Ebenfelt and Rothschild found optimal nondegeneracy conditions
on $M$ and $M'$ such that a germ at $x$ of a local biholomorphism $f$
(between some neighborhoods in $\C^n$)
with ${f(M)\subset M'}$, $f(x)=x'$,
is uniquely determined by its $k$-jet $j^k_xf$,
where $k$ is an integer which depends only on $M$ and $M'$
(see \cite{BER1}, Theorem~1, Proposition~2.3).
A similar statement for the case of hypersurfaces $M$ and $M'$ 
is implicitly contained in \cite{DF1}.
See also \cite{Bel} for the case of Levi-nondegenerate higher
codimensional CR manifolds.
These results should be compared with the following
 theorem of H.~Cartan (\cite{C,N}): \\

{\em Let $D\subset\C^n$ be a bounded domain.
Then the group of all biholomorphic automorphisms $\Aut(D)$ 
equipped with the compact-open topology is a Lie group.
Moreover, if $x\in D$ is arbitrary, an automorphism $f\in \Aut(D)$ is uniquely
determined by its $1$-jet $j^1_xf$ and depends analytically on it.}\\

In the above {\rm CR} setting it is not clear how a map (germ) $f$ depends
on its $k$-jet $j_x^kf$. Is it somehow continuous, analytic etc.?

The special case when both $M$ and $M'$ are Levi-nondegenerate hypersurfaces
was previously considered by Tanaka~\cite{Ta} and Chern-Moser~\cite{CM},
where the determinacy of automorphisms by their $2$-jets was shown.
Even in this special context  
an analytic dependence on the jets is of interest.
An algebraic dependence of the automorphisms of bounded domains 
on their $1$-jets is studied in~\cite{Z},
which can be seen as an algebraic version of the above theorem of Cartan.

The present paper has a different point of view.
Instead of global automorphisms we consider locally
defined germs of holomorphic maps $f$ sending $M$ into $M'$ with no assumptions
on their domains of definition. If either $M$ or $M'$
is locally biholomorphic to $\R^s\times\C^l$, 
the space of such germs is infinite-dimensional 
(see also \cite{BER1}, Theorem~3)
and hence it cannot be parametrized by a 
finite dimensional $k$-jet space.

Therefore both $M$ and $M'$ have to satisfy certain
nondegeneracy conditions. We start with some 
stronger conditions formulated in terms of the Levi form.
Recall that in the case of arbitrary {\rm CR} codimension
the Levi form is a sesqui-linear map (see \S \ref{levi}):
\begin{equation}\label{levi-map}
L\colon T^{1,0}_x M\times T^{1,0}_x M \to (T_x M/T^c_x M)\otimes_{\R}\C.
\end{equation}

\begin{Def}
We say that the Levi form at $x$ is {\em nondegenerate} if for $X_0\in T^{1,0}_xM$,
\begin{equation}
\left( L(X_0,Y)=0 \,\,\, \forall\, Y\in T^{1,0}_xM \right)
\quad \Longrightarrow \quad X_0=0.
\end{equation}
The Levi form is said to be {\em surjective} at $x$,
if the map {\rm(\ref{levi-map})} is surjective.
The Levi form is said to be nondegenerate (resp. surjective)
if it is nondegenerate (resp. surjective) for all $x\in M$.
\end{Def}

Now we introduce the notion of {\em analytic dependency on $k$-jets}.
Let $S$ be a subset of the set of germs at $x\in\C^n$
of all local holomorphic maps $f$ between $\C^n$ and $\C^{n'}$.
By $J^k_x(\C^n,\C^{n'})$ denote the $k$-jet space at $x$
of such maps and by $j^k_xf\in J^k_x(\C^n,\C^{n'})$ the $k$-jet of $f$ at $x$.
We write also $J^k(\C^n,\C^{n'})$ for the (trivial)
bundle of all $k$-jets at all $x\in\C^n$.

\begin{Def}\label{andep}
We say that {\em the germs in $S$ 
depend analytically on their $k$-jets at $x$}
if the following conditions are satisfied:
\begin{enumerate}
\item 
$f\in S$ is uniquely determined within $S$ by $j^k_xf$,
i.e. for $f_1,f_2\in S$, 
$j^k_xf_1=j^k_xf_2 \quad \Longrightarrow \quad f_1=f_2$;
\item 
 For every $f_0\in S$, there exist neighborhoods
$U(x)\subset\C^n$, 
$U(j^k_xf_0)\subset J^k_x(\C^n,\C^{n'})$
such that every $f\in S$ with $j^k_xf\in U(j^k_xf_0)$
extends holomorphically to $U(x)$;
\item In addition there exists 
a holomorphic map
$F\colon U(j^k_xf_0)\times U(x)\to \C^{n'}$
such that for all $f\in S$ with $j^k_xf\in U(j^k_xf_0)$,
\begin{equation}\label{loc-fam}
f(z)=F(j^k_xf,z),\, z\in U(x).
\end{equation}
\end{enumerate}
\end{Def}


For $x\in\C^n$, we equip $S$ with the inductive limit topology, 
i.e. a sequence $(f_n)$ in $S$ converges to $f_0\in S$ 
if and only if all $f_n$ extend to some
neighborhood of $x$ and converge there to $f_0$ uniformly.

\begin{Th}\label{cor1}
        Suppose that the Levi form of $M$ is nondegenerate and surjective.
Then there exists an integer $k>0$ such that for all $x\in M$, 
the germs at $x$ of local biholomorphisms
$f$ with ${f(M)\subset M}$ depend analytically on their $k$-jets.
One can take $k=2(1+\crc M)$.
\end{Th}

As an application we obtain a Lie group structure 
on the set of germs fixing a point $x\in M$ 
(see the end of the section for the proof).

\begin{Cor}\label{cor2}
        Let $x\in M$ be arbitrary and the Levi form of $M$ be nondegenerate
and surjective at $x$. Then the group of all germs of local biholomorphisms 
$f$ with $f(M)\subset M$, $f(x)=x$, is a Lie group.
\end{Cor}

Notice that the domain of definition of a germ $f$ can vary as $j^k_xf$ changes.
It is not even clear a priori whether
a $1$-parameter family $f_t,t\in\R$, of such
automorphisms yields a germ of a vector field. 
The analytic dependence on the $k$-jets
guarantees, in particular,
the following extension result for the germs $f$ with $j_x^kf$ 
sufficiently close to $j_x^kf_0$ for $f_0$ given.

\begin{Cor}
        Under the assumptions of {\rm Theorem~\ref{cor1}}
suppose that $f_0$ is a germ of a local biholomorphism at $x$ with ${f_0(M)\subset M}$.
Then there exists a neighborhood $U(x)\in\C^n$ such that all germs of local biholomorphisms
$f$ with ${f(M)\subset M}$ 
and $j^k_xf$ sufficiently close to $j^k_xf_0$ extend biholomorphically to $U(x)$.
\end{Cor}

Following Tanaka~\cite{Ta} we call a map $\phi\colon M\to M'$
{\em pseudo-conformal}, if it extends to a holomorphic map between
some neighborhoods of $M$ and $M'$ respectively.
One obtains the following global version of Corollary~\ref{cor2}
(see the end of the section for the proof).

\begin{Cor}\label{glob}
        Let $M$ be a {\em compact} {\rm CR} submanifold of 
a complex manifold $X$ which in addition
satisfies the assumptions of {\rm Theorem~\ref{cor1}}.
Then the group of all pseudo-conformal
automorphisms of $M$ is a Lie group.
\end{Cor}

In fact, all of these corollaries 
are proved here in the more general situation,
where $M\subset\C^n$ and $M'\subset\C^{n'}$ 
are generic real-analytic {\rm CR} manifolds 
of arbitrary {\rm CR} dimension and codimension.
In particular, the unique determinacy by $k$-jets, i.e.
the injectivity of the $k$-jet evaluation $f\mapsto j^k_xf$, is also shown.

In this paper $f$ always denotes a germ of a holomorphic map between open
subsets of $\C^n$ and $\C^{n'}$ respectively.
The condition ``$f$ is biholomorphic''
(see \cite{BER},\cite{BER1}) is relaxed to
\begin{equation}\label{cr-surj}
        df_x(T^c_xM) = T^c_{x'}M',
\end{equation}
where $f(x)=x'$ and $T^c_xM := T_xM\cap iT_xM \subset T_x\C^n$
is the complex tangent space.

The following simple example shows
that even if both {\rm CR} dimensions and codimensions of $M$ and $M'$
are equal, the case when (\ref{cr-surj}) 
is satisfied but $f$ is not biholomorphic
is also of interest.\\

\begin{Example}
Let $\phi$ be strongly plurisubharmonic function in a neighborhood of $0\in\C^3$
satisfying
\begin{equation}\label{}
\phi(0)=0, \quad \frac{\partial\phi}{\partial z_3}(0) \ne 0.
\end{equation}
Define $M$ and $M'$ in $\C^4$:
\begin{gather}\label{}
M := \{(z_1,z_2,z_3,z_4): \phi(z_1,z_2,z_3)=\phi(2z_1,z_2,z_4)=0\}, \\
M' := \{(z_1,z_2,z_3,z_4): \phi(z_1,z_2,z_3)= {\rm Re}\, z_4=0\}.
\end{gather}
One has
\begin{equation}\label{}
\crd M = \crd M'=2, \quad \crc M = \crc M'=2.
\end{equation}
The ($\C^2$-valued) Levi forms are given by the second order derivatives
${\partial^2\phi}/{\partial z_i \partial\bar z_j}$.
Simple calculation shows that the Levi form of $M$ (resp. $M'$)
is surjective (resp. nondegenerate).
By Proposition~{\rm \ref{nondeg}} and Lemma~{\rm\ref{minim}},
$M$ and $M'$ satisfy conditions of Theorem~{\rm\ref{main}} below.
The map $f$ given by
\begin{equation}\label{}
f(z_1,z_2,z_3,z_4) := (z_1,z_2,z_3,0),
\end{equation}
satisfies $f(M)\subset M'$ and condition (\ref{cr-surj}) but
is not a diffeomorphism between $M$ and $M'$.
\end{Example}

On the other hand, the condition (\ref{cr-surj}) cannot be removed
as the following elementary example shows.

\begin{Example} 
Define
\begin{equation}
M := \{ |z_1|^2+|z_2|^2-|z_3|^2 = 1 \} \subset \C^3.
\end{equation}
Then $M$ contains the complex line
\begin{equation}
L = \{ z_1 = 1, z_2 = z_3 \}.
\end{equation}
Therefore every holomorphic map 
$f\colon\C^3\to L$
satisfies  $f(M)\subset M$.
Although the Levi form of $M$ is nondegenerate
and surjective, the $k$-jet evaluation $f\mapsto j^k_xf$ is not 
injective for no integer $k$.
\end{Example}

As noted above, the nondegeneracy and surjectivity of the Levi form
are only sufficient conditions.
The optimal conditions which are necessary in many cases
are given in \cite{BER1} for the case $f$ is biholomorphic.
Here we reformulate them in a form suitable for our purposes.

We first recall the notion of the Segre varieties associated
to a generic real-analytic {\rm CR} submanifold $M\subset\C^n$.
Let $x\in M$ be arbitrary and let $M$ be defined near $x$ 
by the real-analytic equations

\begin{equation}\label{13}
\rho_1(z,\bar z) = \cdots = \rho_d(z,\bar z) = 0, 
\quad \partial\rho_1\wedge\cdots\wedge\partial\rho_d \ne 0.
\end{equation}

The complexification $\M\subset\C^n\times\ov{\C^n}$ is defined by
\begin{equation}\label{14}
\M := \{ (z,\bar w) \in U(x)\times U(\bar x) : 
\rho_1(z,\bar w) = \cdots = \rho_d(z,\bar w) = 0 \},
\end{equation}
where $U(x)\subset\C^n$ and $U(\bar x)\subset\ov{\C^n}$
are sufficiently small neighborhoods.
As a germ at $(x,\bar x)$ of a complex-analytic subset,
$\M$ is independent of the choice of the defining equations (\ref{13}).
By fixing the coordinate $\bar w$, we obtain the Segre variety
\begin{equation}\label{15}
Q_w := \{ z\in U(x) : (z, \bar w) \in \M \}.
\end{equation}
Segre varieties were introduced by Segre \cite{Se} and play an important
role in the reflection principle (see e.g. \cite{P,L,DW,W,W1}).
For $z,w$ close to $x$, it follows that $Q_w$ is a complex manifold of dimension $n-d=\crd M$.
The following symmetry property is a direct corollary of the invariance
of $\M$ under the involution $(z,\bar w)\mapsto (w,\bar z)$:
\begin{equation}
z\in Q_w \iff w\in Q_z.
\end{equation}

For arbitrary $x\in\C^n$,
denote by $J^{k,d}_x(\C^n)$ the space of $k$-jets at $x$
of $d$-codimensional complex submanifolds $V\subset\C^n$ with $x\in V$
and by $j^k_x(V)\in J^{k,d}_x(\C^n)$ the $k$-jet of $V$.
The following is another form of the definition of a $k$-nondegenerate
manifold given in \cite{BER}.

\begin{Def}\label{knondeg}
        A generic real-analytic {\rm CR} manifold $M\subset\C^n$ is called
{\em $k$-nondegenerate} at $x\in M$ ($k\ge 1$), if the (antiholomorphic) map
\begin{equation}
\phi^k\colon Q_x\to J^{k,d}_x(\C^n),\, w\mapsto j^k_x(Q_w),
\end{equation}
is of the rank $n-d = \dim Q_x$ at $x$. $M$ is called {\em finitely nondegenerate},
if it is $k$-nondegenerate for some $k$.
\end{Def}

A manifold $M$ is called {\em essentially finite} at $x$
if the germ of $Q_x$ at $x$ is different from germs of $Q_w$ at $x$
for all $w\in Q_x\setminus\{x\}$ sufficiently close to $x$.
This notion was introduced in \cite{BJT} and earlier implicitly
discussed in \cite{DF1}.
By Proposition~1.3.1 in \cite{BER}, if $M$ is essentially finite at $x$,
there exists the so-called Levi number $l(M)$. Then
$M$ is $l(M)$-nondegenerate at all points $p\in M\setminus W$,
where $W\subset M$ is a proper real-analytic subset.
Clearly the $k$-nondegeneracy guarantees
the $k'$-nondegeneracy for all $k'\ge k$.

The second important notion is the minimality condition
introduced by Tumanov in \cite{Tu}:

\begin{Def}\label{min}
        A {\rm CR} manifold $M$ is called {\em minimal } at $x\in M$,
if there does not exist a proper {\rm CR} submanifold $N\subset M$ with $x\in N$
such that $\crd N = \crd M$.
\end{Def}

Tumanov \cite{Tu} shows that, if $M$ is minimal at $x$,
all {\rm CR} functions on $M$ extend holomorphically to a wedge with the edge $M$.

\begin{Def}\label{admissible}
We call a germ $f$ {\em admissible}
if it satisfies the following conditions:
\begin{equation}\label{admissibleeq}
f(M)\subset M', \quad df(T^c_x M) = T^c_{f(x)} M'.
\end{equation}
\end{Def}

\begin{Th}\label{main}
        Suppose that $M'$ is $r$-nondegenerate,
$M$ is minimal at $x$ and $k:=2r(1+\crc M)$.
Let $S$ be the set of all admissible germs $f$ at $x$.
Then the germs in $S$ depend analytically on their $k$-jets at $x$.\par
\end{Th}

A different proof of Theorem~\ref{main} in the case of hypersurfaces 
($\crc M = \crc M' = 1$)
with a sharper estimate for the jet order ($2r$-jets instead of $4r$-jets)
has been obtained recently by Baouendi, Ebenfelt and Rothschild \cite{BER2}.\\

\noindent {\bf Remark.}
The case, where one of the manifolds 
$M,M'$ is not generic, can be reduced to the
generic case. Every real-analytic {\rm CR} submanifold $M\subset\C^n$ 
is generic in
the so-called {\em intrinsic complexification} $V\subset\C^n$ defined to be
the minimal (in the sense of germs)
 complex-analytic subvariety which contains $M$. 
If $M$ is \CR, $V$ is smooth. After a change of
local coordinates $V$ becomes open in $\C^m\subset\C^n$ near $x\in M$.
The condition $f(M)\subset M'$ automatically implies $f(V)\subset V'$.
If $M'$ is not generic, we can therefore replace $\C^{n'}$ with $V'$.
If $M$ is not generic, Theorem~\ref{main} yields 
a parametrization (\ref{loc-fam})
for the restricted admissible germs $f\colon V\to\C^{n'}$,
whereas outside $V$ the germs can be chosen arbitrary
and cannot be determined by their $s$-jets even for $s$ arbitrary large.

\vspace{5pt}

\noindent As a corollary, we obtain the following description
of the space of admissible germs as a (locally closed) real-analytic subset of
$J^k_x(\C^n,\C^{n'})$.

\begin{Cor}\label{cor-main}
        Let $k,r,x,M,M',S$ satisfy the conditions of {\rm Theorem~\ref{main}}.
\begin{enumerate}
\item Then the set $A\subset J^k_x(\C^n,\C^{n'})$
of the $k$-jets of all $f\in S$ is a locally closed real-analytic subset.
\item 
There exists a neighborhood $U$ of
$A\times \{x\} \subset A\times \C^n$
and a (unique) real-analytic map
\begin{equation}\label{a-fam}
        G \colon U \to \C^{n'},
\end{equation}
which is holomorphic in the $\C^n$-factor and is
such that every $f\in S$ satisfies
\begin{equation}\label{f-fam}
        f(w) =  G (j^k_x f, w).
\end{equation}
\end{enumerate}
\end{Cor}

\begin{proof}
        By Theorem~\ref{main}, the germs in $S$ depend analytically
on their $k$-jets.
Let $M,M'$ be closed in some open subsets 
\begin{equation}
U(M)\subset\C^n,\quad U(M')\subset\C^{n'}
\end{equation}
respectively. 
For every $f_0\in S$, Theorem~\ref{main} yields the local
parametrization (\ref{loc-fam}). 
Without loss of generality, $U(x)\subset U(M)$ and 
\begin{equation}
F(U(j^k_xf_0)\times U(x))\subset U(M').
\end{equation}
Then the above set $A\subset J^k_x(\C^n,\C^{n'})$ is locally defined by
\begin{multline}
A\cap U(j^k_xf_0) := \{ j\in U(j^k_xf_0) : \\ F(j,w)\subset M'
\text{ for all } w\in M\cap U(x), \quad j=j^k_x F(j,\cdot) \}.
\end{multline}
It follows from the real-analyticity of $M$ and $M'$ that $A$ is locally
defined by finitely many real-analytic equations,
i.e. it is a real-analytic subset.
The restrictions of $F$'s to $A\times U(x)$ for different $f_0$'s
glue together to a well-defined real-analytic map $G$ 
satisfying the required properties.
\end{proof}

\begin{proof}[Proof of Corollary~\ref{cor2}]
By Corollary~\ref{cor-main}, the set $A\subset J^k_x(\C^n,\C^{n'})$ 
of the $k$-jets ($k:={2(d+1)r}$) of all admissible germs is real-analytic.
Denote by $A_0\subset A$ the subset of the $k$-jets fixing $x$. 
Clearly $A_0$ is also real-analytic and locally closed.
Since $A_0$ is also closed under the composition of the jets,
it is a closed Lie subgroup of the Lie group of all $k$-jets fixing $x$.
\end{proof}

\begin{proof}[Proof of Corollary~\ref{glob}]
By Theorem~\ref{main} and the definition 
of analytic dependence on $k$-jets,
there exists a finite open covering of $M$ with
coordinate neighborhoods $U(x_i)$ with parametrizations
\begin{equation}
F_i\colon U(j^k_{x_i}f_0)\times U(x_i)\to \C^n,
\end{equation}
where $f_0:=\id\in \Aut(M)$.
An automorphism $f\in \Aut(M)$ in a neighborhood of $\id$ is 
therefore uniquely
determined by finitely many parameters
\begin{equation}\label{}
p_i=j^k_{x_i}f, \quad q_i= j^k_{x_i}f^{-1}.
\end{equation}
Set $f_i:=F_i(p_i,\cdot)$, $g_i:=F_i(q_i,\cdot)$.
A tuple $(p_1,\ldots,p_m,q_1,\ldots,q_m)$ determines an automorphism
of $M$ if and only if the following is satisfied:
\begin{enumerate}
\item $f_i=f_j$, $g_i=g_j$ on $U(x_i)\cap U(x_j)$ for all $i,j$;
\item $f_i(M)\subset M$, $g_i(M)\subset M$ for all $i$;
\item The glued maps $f_p$ and $g_q$ satisfy
        $f_p\circ g_q = g_q\circ f_p = \id$, $j^k_{x_i}f_p=p_i$, $j^k_{x_i}g_q=q_i$.
\end{enumerate}
Conditions (1.),(2.),(3.) define an analytic subset of the parameter space
of $p_i$'s and $q_i$'s.
In these coordinates the group operation is given by 
$$ ((p,q),(p',q')) \mapsto (j^k_{x_i}(f_{p'}\circ f_p),j^k_{x_i}(g_q\circ g_{q'}))$$
and is therefore real-analytic. Hence $\Aut(M)$ is a Lie group.
\end{proof}

\section{The Levi form and nondegeneracy conditions.}\label{levi}

Let $M\subset\C^n$ be a CR submanifold and $x\in M$ an arbitrary point.
Recall that a $(1,0)$-vector field on $M$ is a vector field $X$ in
${T^cM\otimes_{\R}\C}$ such that $JX = iX$,
where $J\colon T^cM\otimes_{\R}\C \to T^cM\otimes_{\R}\C$
is the complexification of the CR structure
$J\colon T^cM\to T^cM$ and $iX$ is the multiplication by $i$
in the component $\C$ of the tensor product.

The Levi form of $M$ at $x$ is the hermitian (vector-valued) form
\begin{equation}
L\colon T^{1,0}_xM\times T^{1,0}_xM \to (T_xM/T^c_xM)\otimes_{\R}\C,
\end{equation}
given by
\begin{equation}
L(X,Y) := \frac{1}{2i} \pi([X,\bar Y]),
\end{equation}
where $X,Y$ are $(1,0)$-vector fields on $M$ and
\begin{equation}
\pi\colon T_xM\otimes_{\R}\C \to (T_xM/T^c_xM)\otimes_{\R}\C
\end{equation}
is the canonical projection.
Notice that $L(X,Y)(x)$ depends only on $X(x),Y(x)$.

\smallskip
\noindent {\bf Remark.}
If $M\subset\C^n$ is a hypersurface defined by
$M=\{\phi=0\}$ with $d\phi\ne 0$,
the standard Levi form of $\phi$ coincides 
with the evaluation of $\partial\phi$ on $L(X,Y)$.
\smallskip

The Levi form is a first order holomorphic invariant of $M$
whereas the conditions in Theorem~\ref{main} are of possibly higher order.
However, we obtain the following

\begin{Prop}\label{nondeg}
        A generic real-analytic ${\rm CR}$ manifold $M\subset\C^n$
is $1$-nondegenerate at $x$
if and only if the Levi form at $x$ is nondegenerate.
\end{Prop}

\begin{proof}
        Let $M$ be locally defined near $x$ by the real-analytic equations
\begin{equation}\label{rho}
\rho_1(z,\bar z) = \cdots = \rho_d(z,\bar z) = 0,
\end{equation}
where $\partial\rho_1 \wedge \cdots \wedge \partial\rho_d \ne 0$.

The differentials $\partial\rho_j$ vanish on $T^cM\otimes_{\R}\C$.
Hence, for a $(1,0)$-vector field $X\subset T^{1,0}M$,
\begin{equation}
\partial\rho_j(X) = \partial\rho_j(\bar X) = 0.
\end{equation}
Denote by $X\psi$ the derivative of $\psi$ along $X$.

Applying the formula
\begin{equation}\label{}
d\omega(X_1,X_2) = X_1 \omega(X_2) - X_2 \omega(X_1) - \omega([X_1,X_2])
\end{equation}
for $\omega = \partial\rho_j$, $X_1=X, X_2=\bar Y$, we obtain
\begin{equation}\label{levi-con}
\bar\partial\partial\rho_j (X,\bar Y) = - \partial\rho_j([X,\bar Y]).
\end{equation}

The Segre varieties $Q_w$ are given by
\begin{equation}
\rho_1(z,\bar w) = \cdots = \rho_d(z,\bar w) = 0
\end{equation}
and their $1$-jets $j^1_xQ_w$ by the linear equations
\begin{equation}
\partial\rho_1(z,\bar w) = \cdots = \partial\rho_s(z,\bar w) = 0.
\end{equation}

Let $e_1(z,\bar w),\ldots,e_d(z,\bar w)$ be a collection of $(1,0)$-vector 
fields on the complexification $\M$ which is
a basis of $T_zQ_w$ at every point $(z,\bar w)\in \M$ close to
$(x,\bar x)$. Let $e'_1,\ldots,e'_d$ be a similar collection
of $(1,0)$-vector fields in the $w$-direction, i.e. pointwise a basis of
$T_{\bar w}\ov{Q_{\bar z}}$.
The rank of the map $\phi^1$ in Definition~\ref{knondeg} is the same as
the rank of the matrix with the rows
\begin{equation}
\bar\partial\partial\rho_j (e_k,\bar e'_l) \text{ with $l$ fixed }.
\end{equation}

By (\ref{levi-con}), this matrix has the rows
\begin{equation}\label{matr}
-\partial\rho_j ([e_k,\bar e'_l]) \text{ with $l$ fixed }.
\end{equation}

Since $\partial\rho_1 \wedge \cdots \wedge \partial\rho_d \ne 0$,
the tuple $(\partial\rho_1,\ldots,\partial\rho_d)$ defines an isomorphism
between $(TM/T^cM)\otimes_{\R}\C$ and $\C^d$.
Hence the matrix with the rows (\ref{matr}) has the same rank as
the matrix with the rows
\begin{equation}
L(e_k,e_l) \text{ with $l$ fixed }.
\end{equation}
This rank equals to $n-d$
if and only if the Levi form of $M$ is nondegenerate.
\end{proof}

The minimality condition (Definition~\ref{min})
involves high order commutators 
and therefore cannot be formulated in terms of the Levi form.
However, one obviously has the following sufficient condition.

\begin{Lemma}\label{minim}
        Suppose that the Levi form of $M$ at $x$ is surjective onto
$(T_xM/T^c_xM)\otimes_{\R}\C$. Then $M$ is minimal at $x$.
\end{Lemma}

Now we can formulate special cases of Theorem~\ref{main}
and Corollary~\ref{cor-main} under Levi form conditions.

\begin{Th}\label{levi-thm}
        Suppose that the Levi form of $M'$ is everywhere nondegenerate,
the Levi form of $M$ is surjective at $x$ and $d=\crc M$.
Then the conclusions of {\rm Theorem~\ref{main}}
 and {\rm Corollary~\ref{cor-main}} hold.
\end{Th}

\section{Local parametrization of jets of holomorphic maps.}

One of the technical tools for proving Theorem~\ref{main}
is the following connection between the $k$-jets at $z$
and the $(k+r)$-jets at $\bar w$.
For simplicity, we write for a subset $U(j)\subset J^{k+r}(\C^n,\C^{n'})$,
\begin{equation}
  U_z(j) := U(j) \cap J_z^{k+r}(\C^n,\C^{n'}).
\end{equation}

\begin{Prop}\label{fam}
       Suppose that $M'$ is $r$-nondegenerate at $x'$ and
$f_0$ is an admissible germ at $x\in M$, $x'=f_0(x)$.
Then, for every integer $k$ there exist neighborhoods
$U(j_x^{k+r}f_0)\subset J^{k+r}(\C^n,\C^{n'})$, 
$U(x,\bar x)\subset \M$
and for every $(z,\bar w)\in U(x,\bar x)$ a holomorphic map
        $F_{(z,\bar w)}^k \colon U_z(j_x^{k+r}f_0) \to
        J^k_{\bar w}(\ov{\C^n},\ov{\C^{n'}})$,
such that for all admissible germs $f$ at $x$ with
$j^{k+r}_x f\in U(j^{k+r}_x f_0)$, 
\begin{equation}\label{sat}
        j^k_{\bar w} \bar f =  F_{(z,\bar w)}^k (j^{k+r}_z f),
\end{equation}
where $(z,\bar w)\subset\M$ is sufficiently close to $(x,\bar x)$.
Moreover, the map $F_{(z,\bar w)}^k$
depends holomorphically on $(z,\bar w)\in U(x,\bar x)$.
\end{Prop}

A different proof of Proposition~\ref{fam} is given in \cite{BER} 
(Assertion~3.3.1) and \cite{BER1} (Proposition~2.2).

\begin{proof}
Let $(z,\bar w)\in\M$ be close to $(x,\bar x)$.
By the construction (see (\ref{13}),(\ref{14}),(\ref{15})),
the Segre variety $Q_w$ is smooth at $z$
and $d:=\dim_z Q_w$ is constant. 
We choose local coordinates 
$z=(z_1,z_2)\in\C^d\times\C^{n-d}$
near $x$ such that
the Segre variety $Q_w$ has the form of the graph:
\begin{equation}\label{graph}
        Q_w = \{ (z_1,\phi_w(z_1)) : z_1\in U(x_1) \}, \quad
        x = (x_1,x_2),
\end{equation}
where
$U(x_1)\subset\C^n$ is an open neighborhood and
        $\phi_w \colon U(x_1)\to \C^{n-d}$
is holomorphic.
The map $\phi_w$ is uniquely defined and depends
holomorphically on $\bar w$ in some neighborhood 
$U(\bar x)\subset\ov{\C^n}$ .
Denote by $j^r_z Q$ the holomorphic $r$-jet evaluation
\begin{equation}
        j^r_z Q \colon U(\bar x) \to J^r_{z_1}(\C^d,\C^{n-d}), \quad
        \bar w \mapsto j^r_{z_1} \phi_w.
\end{equation}

A similar evaluation is obtained for $(z',\bar w')\in \M'$
sufficiently close to $(x',\bar x')$:
\begin{equation}\label{Q'-map}
        j^r_{z'} Q' \colon U(\bar x') \to J^r_{z'_1}(\C^{d'},\C^{n'-d'}),\quad
        \bar w' \mapsto j^r_{z'_1} \phi'_{w'}.
\end{equation}

We claim that (\ref{Q'-map}) is an immersion.
Since $M'$ is $r$-nondegenerate at $x'=(x'_1,x'_2)$,
$j^r_x Q'$ restricted to the Segre variety $Q'_{x'}$ is an immersion.
By the construction, the $0$-jet evaluation 
$j^0_{x'} Q' = x'_2$ is constant on $Q'_{x'}$.
On the other hand, the restriction of $j^0_{x'} Q'$ to 
the transversal direction, i.e. $\{x'_1\}\times\C^{n'-d'}$,
is also an immersion.
Indeed, since
$M'\cap (\{x'_1\}\times\C^{n'-d'})$
is totally real in $\{x'_1\}\times \C^{n'-d'}$, 
it is locally biholomorphically equivalent to 
$\R^{n'-d'}$.
For $\R^s\subset\C^s$, 
the immersion property of $j_{x'}^0Q'$ can be directly verified.

Hence the total map (\ref{Q'-map}) splits into
$j^0_{x'}Q'$ which is an immersion on $\{x'_1\}\times \C^{n'-d'}$
and constant on $Q'_{x'}$
and the remainder which is an immersion on 
$Q'_{x'}$. Therefore (\ref{Q'-map}) is an immersion
for $z$ close to $x$ and 
locally there exists a left inverse map $(j^r_{z'} Q')^{-1}$ which satisfies
   $(j^r_{z'} Q')^{-1} \circ (j^r_{z'} Q') = {\rm \id}$.

Let $f_0$ be an admissible germ.
Complexifying the condition $f_0(M)\subset M'$ we obtain
$(f_0,\bar f_0)(\M) \subset \M'$,
which means
\begin{equation}\label{inc}
        f_0(Q_w) \subset Q'_{f_0(w)}.
\end{equation}

Since $\dim Q_w$ is constant, the second condition $df(T^c_xM)=T^c_{f(x)}M'$
implies
\begin{multline}
\dim f_0(Q_w) = \dim df_0(T_xQ_x) = \dim df_0(T^cM) \\ =
\dim T^cM' = \dim Q_{x'} = \dim Q'_{f_0(w)}.
\end{multline}
Along with (\ref{inc}) this yields
\begin{equation}\label{51}
        f_0(Q_w) = Q'_{f_0(w)}, \quad df_0(T_zQ_w) = T_{z'}Q'_{w'}. 
\end{equation}

After going to the $k$-jet evaluations 
we obtain the following commutative diagram
\begin{equation}\label{52}
\begin{CD}
U(\bar x)                       @>{j^r_z Q}>>   J^r_{z_1}(\C^d,U(x_2)) \\
@V{\bar f_0}VV                                    @VV{f_{0*}}V            \\
U(\bar x')              @>{j^r_{z'} Q'}>>    J^r_{z'_1}(\C^{d'},U(x'_2)),
\end{CD}
\end{equation}
where $f_{0*}$ is the corresponding map on the jet level defined as follows.

By (\ref{51}), there exists a local splitting of the source space
$\C^d = \C^{d'} \times \C^{d-d'}$
near $x_1=(x_{11},x_{12})$ such that the restriction
\begin{equation}\label{1}
f_0 \colon Q_w\cap \{z_{12}=x_{12}\} \to Q'_{w'}
\end{equation}
is locally biholomorphic for all $(z,\bar w)\in \M$
close to $(x,\bar x)$.

We write
$f_{01}\colon U(x) \to \C^{d'}$, 
$f_{02}\colon U(x) \to \C^{n'-d'}$
for the components of $f_0$ 
(here $\C^{d'}$ and $\C^{n'-d'}$ are in the target space).
Since (\ref{1}) is locally biholomorphic, the map
\begin{equation}\label{2}
\tilde f_0 \colon \C^{d'}\to \C^{d'}, \quad
z_{11} \mapsto f_{01}(z_{11},x_{12},\phi_{\bar w}(z_{11},x_{12}))
\end{equation}
is also locally biholomorphic.

Let $\tilde f_0^{-1}$ be the local inverse.
By (\ref{graph}) and (\ref{51})
\begin{equation}\label{2.1}
\phi'_{\bar w'}(z'_1) = 
f_{02}(\tilde f_0^{-1}(z'_1),x_{12},
\phi_{\bar w}(\tilde f_0^{-1}(z'_1),x_{12}))
\end{equation}
and hence, passing to the $r$-jets,
\begin{equation}\label{3}
j^r_{z'_1} \phi'_{\bar w'} = 
(j^r_z f_{02})( (j^r_{z_{11}}\tilde f_0)^{-1},x_{12},
(j^r_{z_1}\phi_{\bar w})((j^r_{z_{11}}\tilde f_0)^{-1},x_{12})),
\end{equation}
where by (\ref{2}),
\begin{equation}\label{4}
j^r_{z_{11}}\tilde f_0 =
(j^r_z f_{01}) (z_{11},x_{12},
(j^r_{z_1}\phi_{\bar w}) (z_{11},x_{12}))
\end{equation}
with obvious notation.

The formulae (\ref{3},\ref{4}) can be written together in the form
\begin{equation}\label{5}
j^r_{z'_1} \phi'_{\bar w'} = 
\Phi (j^r_z f_0, j^r_{z_1}\phi_{\bar w}),
\end{equation}  
where
\begin{equation}
\Phi \colon U_z(j^r_x f_0) 
\times U_{z_1}(j^r_{x_1}\phi_{\bar x})
\to J^r_{z'_1}(\C^{d'},U(x'_2))
\end{equation} 
is a holomorphic family of maps and
$U(j^r_x f_0)\subset J^r(\C^n,\C^{n'})$,
$U(j^r_{x_1}\phi_{\bar x})\subset J^r(\C^d,\C^{n-d})$
are sufficiently small neighborhoods.

Define $f_{0*}$ as in (\ref{52}) by
\begin{equation}\label{6}
f_{0*} := \Phi (j^r_z f_0, \cdot).
\end{equation}
By (\ref{5}), the diagram (\ref{52}) is commutative.
\newline
The com\-mu\-ta\-ti\-vi\-ty means that
\begin{equation}
(j^r_{z'} Q') \circ \bar f_0 = f_{0*} \circ (j^r_z Q).
\end{equation}
Applying the left inverse $(j^r_{z'} Q')^{-1}$ to both sides we obtain
\begin{equation}\label{comp}
\bar f_0 = (j^r_{z'} Q')^{-1} \circ f_{0*} \circ (j^r_z Q)
\end{equation}
and, passing to the $k$-jets,
\begin{equation}\label{7}
j^k_{\bar w} \bar f_0 = 
j^k_{g'}(j^r_{z'} Q')^{-1} \circ j^k_g f_{0*} \circ 
j^k_{\bar w} (j^r_z Q),
\end{equation}
where $g :=  (j^r_z Q) (\bar w)$ and $g' := j^k_g f_{0*} (g)$.

The $k$-jet $j^k_g f_{0*}$ can be calculated from (\ref{6}):
\begin{equation}\label{8}
j^k_g f_{0*} := (j^k_s \Phi) (j^{k+r}_z f_0, \cdot), \quad
s := (j^r_z f_0, g).
\end{equation} 

We write (\ref{7}) and (\ref{8}) together in the form
\begin{equation}\label{9}
j^k_{\bar w} \bar f_0 = F^k_{z,\bar w} (j^{k+r}_z f_0),
\end{equation} 
where 
\begin{equation}
 F^k_{z,\bar w} \colon U_z(j^{k+r}_x f_0) \to 
J^k_{\bar w}(\ov{\C^n},\ov{\C^{n'}})
\end{equation} 
is a holomorphic family of maps and
$U(j^{k+r}_x f_0)\subset J^{k+r}(\C^n,\C^{n'})$
a sufficiently small open neighborhood.

Now we notice that $f_0$ was an arbitrary admissible germ.
Let $f$ be another one with 
$j^{k+r}_x f \in U(j^{k+r}_x f_0)$.
Then for $(z,\bar w)\in \M$ sufficiently close to $(x,\bar x)$,
(\ref{2.1})-(\ref{9}) remain valid if we replace $f_0$ with $f$.
In particular, the proposition follows from (\ref{9})
with $f_0$ replaced by $f$.
\end{proof}

\section{Proof of Theorem~\ref{main}}

Given $j^{k+r}_z$ for $z\in\C^n$ fixed,
(\ref{sat}) yields the $k$-jet $j^k_wf$, in particular $f$ itself
($f(w)=j^0_wf$). However, the point $w$ is not arbitrary
but must satisfy
\begin{equation}
(z,\bar w)\in\M, \text{ i.e. } w\in Q_z.
\end{equation}

Our next goal 
will be to obtain a similar formula
which yields $f$ at least on an open subset of $\C^n$.
The idea here is to iterate the construction of the previous
paragraph.

\subsection{Iterated complexification.}

The iterated complexification $\M^s$ of $\M$ is defined by
\begin{multline}\label{odd}
        \M^{2l-1} := \{(z_1,\bar w_1,\ldots, z_l,\bar w_l)
        \in \C^n\times\ov{\C^n}\times\cdots\times\C^n\times\ov{\C^n}
        : \\ (z_1,\bar w_1)\in\M, \, (w_1,\bar z_2)\in\M, \,
        \ldots, \, (w_{l-1},\bar z_l)\in\M, \, (z_l,\bar w_l)\in\M \},
\end{multline}
if $s=2l-1$ is odd and by
\begin{multline}\label{even}
        \M^{2l} := \{(z_1,\bar w_1,\ldots, z_l,\bar w_l, z_{l+1})
        \in \C^n\times\ov{\C^n}\times\cdots\times\C^n\times\ov{\C^n}\times\C^n
        : \\ (z_1,\bar w_1)\in\M,\, (w_1,\bar z_2)\in\M,\,
        \ldots,\, (z_l,\bar w_l)\in\M,\, (w_l,\bar z_{l+1})\in\M \},
\end{multline}
if $s=2l$ is even.

For simplicity we restrict the proof of Theorem~\ref{main} in the sequel
to the case (\ref{odd}) of odd iterates. 
The case (\ref{even}) of even iterates
is completely analogous.
We write
\begin{equation}
        \M^s_z := \M^s \cap \{z_1 = z\}.
\end{equation}

\begin{Lemma}\label{M-smooth}
The analytic subset $\M^s\subset\C^N$ is smooth at $(x,\ldots,\bar x)$
and the projection on the first copy of $\C^n$ is of rank $n$ at
this point. $\M^s_x$ is also smooth at $(x,\ldots,\bar x)$.
\end{Lemma}

\begin{proof}
By the construction of the Segre varieties 
(see (\ref{13}),(\ref{14}),(\ref{15})), there exists a local splitting
$\ov{\C^n(w)}=\ov{\C^d(w^1)}\times\ov{\C^{n-d}(w^2)}$
and a germ of a local holomorphic map
$\phi\colon \C^n \times \ov{\C^d} \to \ov{\C^{n-d}}$
such that the germ 
\begin{equation}
\C^n \times \ov{\C^d} \to \M, \quad
(z,\bar w^1)\mapsto (z,\bar w^1, \phi(z,\bar w^1)),
\end{equation}
is biholomorphic. Then the germ
\begin{multline}
\Phi\colon \C^n \times \ov{\C^d} \times \C^d \times \cdots \times \ov{\C^d}  \to \M^s,\\
(z_1,\bar w^1_1,\ldots,z^1_l,\bar w^1_l)\mapsto \\
(z_1,\bar w^1_1, \phi(z_1,\bar w^1_1)), z^1_2, 
\bar\phi(\bar w^1_1, \phi(z_1,\bar w^1_1), z^1_2), \ldots,
\bar w^1_l, \phi(\ldots, \bar w^1_l)),
\end{multline}
is also biholomorphic, where $s=2l-1$. 
Furthermore, the subset $\M^s_x$ is equal to $\Phi(\{z_1=x_1\})$.
This proves the lemma.
\end{proof}

\begin{Lemma}\label{M-lemma}
Under the assumptions of {\rm Proposition~\ref{fam}},
for every integers $s\ge 1$, $l\ge 0$, $m:=sr+l$, 
there exist open neighborhoods
\begin{equation}
U(j_x^mf_0)\subset J^m(\C^n,\C^{n'}), \quad
U(x,\bar x,\ldots,\bar x)\subset \M^s
\end{equation}
 and for every 
$\nu=(z,\bar w_1,\ldots,z_l,\bar w)\in\M^s$,
a family of holomorphic maps
\begin{equation}\label{Fnu}
        F^l_\nu \colon U_z(j_x^mf_0) \to J^l_{\bar w}(\ov{\C^n},\ov{\C^{n'}})
\end{equation}
such that for all admissible germs $f$ at $x$ with
$j^m_x f\in U(j^m_x f_0)$,
\begin{equation}\label{f-F}
        j^l_{\bar w}\bar f = F^l_\nu( j^m_z f ),
\end{equation}
where $z,w$ and all $z_j,w_j$ are sufficiently close to $x$.
Moreover, $F_\nu$ depends holomorphically on $\nu\in\M^s$.
\end{Lemma}

\begin{proof}
The required maps $F^l_\nu, \nu\in\M^s$,
are obtained as iterates of (\ref{sat}):
\begin{equation}\label{sat-it}
F_\nu (j) := \ov {F_{(w_l,\bar z_{l+1})}^l} \circ F_{(z_l,\bar w_l)}^{r+l}
        \circ \cdots \circ F_{(z_1,\bar w_1)}^{(s-1)r+l} (j).
\end{equation}

Here we iterate step by step the map in (\ref{sat}) and its conjugate.
Notice that the conjugate is taken exactly for $z_j$'s and $w_j$'s
conjugated to the coordinates of $\nu$.
Hence $F_\nu$ depends holomorphically on $\nu\in\M^s$.
The required formula (\ref{f-F}) 
for admissible germs  is obtained
by iterating  (\ref{sat}).
\end{proof}

\subsection{Segre sets.}

The right-hand side in (\ref{f-F}) depends
on $\nu$ which contains several coordinates $z_j,w_j$
other than $z$ and $w$.
To avoid this ambiguity we project out the auxiliary coordinates.
The corresponding projection of $\M^s$ is a family of
the so-called 
{\em Segre sets} introduced in \cite{BER}.
The family of Segre sets
\begin{equation}
        Q^{2l+1}\subset\C^n\times\ov{\C^n},\,
        Q^{2l}\subset\C^n\times\C^n
\end{equation}
is defined as the projection
on the product of the fist space $\C^n$ and the last one
(which is either $\ov{\C^n}$ or $\C^n$)
of $\M^s$ for $s=2l+1$ or $s=2l$ respectively.

Notice that, on the contrary to the family of Segre varieties $Q_w = Q^1_w$,
the family of Segre sets $Q^s_w$ for $s>1$ is not analytic in general
and depends on the neighborhood, where $\M$ is taken.
A smaller neighborhood of $(x,\bar x)\in \M$ induces
a family of smaller Segre sets in {\em each} neighborhood of 
$(x,\bar x)$.

Let 
\begin{equation}\label{qsw}
Q^s_w := \{z\in\C^n : (z,\bar w)\in Q^s \}
\end{equation}
be the $s$-th Segre set associated with $w\in\C^n$.
We make use of the following basic properties of the Segre sets
(see \cite{BER}, \S2.2, in particular Corollary~2.2.2):  

\begin{Th}[Baouendi, Ebenfelt, Rothschild]\label{opensegre}
        Let $M\subset\C^n$ be a real-analytic CR submanifold.
Then the dimension of $Q^s_z$ increases strictly until it stabilizes.
If $M$ is minimal at $x\in M$, then $Q^s_x$ contains an open subset
of $\C^n$ for $s$ sufficiently large.
\end{Th}

Since 
\begin{equation}
\dim Q^1_w =\dim Q_w = \crd M,
\end{equation}
we can take $s=d+1$, where $d=\crc M$.
By Theorem~\ref{opensegre}, the generic rank of the projection 
\begin{equation}\label{pr}
\M^{d+1}_x\to \ov{\C^n}
\end{equation}
is $n$.
Unfortunately this is true only generically whereas
the rank of (\ref{pr}) at $(x,\bar x, \ldots,\bar x)$ is usually not maximal.
Hence we cannot lift $Q^s_x$ to $M$ in a neighborhood of $x$
but we shall do this in an open set of $\ov{\C^n}$ with $\bar x$ on the boundary.

\subsection{Liftings of the Segre sets to the iterated complexification.}

Our next goal will be to obtain a formula similar to (\ref{f-F})
where the family of holomorphic maps is parametrized by $w\in Q_z^s$
instead of $\nu\in \M^s$. For this we choose carefully open subsets
of the Segre sets $Q_z^s$ (see (\ref{qsw})) and lift them simultaneously to $\M_z^s$.
The pullbacks of (\ref{Fnu}) under these liftings
will yield the required formulae.

We start with some elementary lemmata which we 
prove here for convenience of the reader.

\begin{Lemma}\label{eq-dim}
Let $X$, $Y$ be complex manifolds, $x_0\in X$, $\dim X = \dim Y = n\ge 1$,
$Y'\subset Y$ a submanifold
and $f\colon X\to Y$ a holomorphic map
with $y_0:=f(x_0)\in Y'$.
Suppose that $f$ is of the maximal rank $n$ on a dense subset of $f^{-1}(Y')$.
Then there exists an open subset $N\subset Y$ and a holomorphic lifting
$\psi\colon N\to X$ (i.e. $f\circ\psi=\id_N$) 
such that the following is satisfied:
\begin{enumerate}
\item $N\cap Y'$ is connected;
\item $y_0$ is in the closure of $N\cap Y'$;
\item $\psi$ extends continuously to $y_0$ with $\psi(y_0)=x_0$.
\end{enumerate} 
\end{Lemma}

\begin{proof}
Without loss of generality, $X$ is an open subset in $\C^n$.
Set
\begin{equation}
D:= \{ x\in X : \rk_xf < n \}.
\end{equation}
Since $D\subset X$ is an analytic subset and $D\cap f^{-1}(Y')$
is nowhere dense, there exists a local holomorphic curve 
$C\subset f^{-1}(Y')$ such that $C\cap C\subset \{y_0\}\subset C$.
Then the restriction of $f$ to $C$ is a ramified covering 
given by $z\mapsto z^k$ with respect to some local coordinates.
Therefore there exists a real curve
\begin{equation}
\gamma\colon [0,1] \to C, \quad \gamma (0) = x_0,
\quad f \text{ is injective on } \gamma[0,1].
\end{equation}

Let $t_0\in[0,1)$ be minimal with the property that there exist 
$N$ and $\psi$ satisfying (1.) and (3.) in this lemma and such that
\begin{equation}
f(\gamma(t))\subset N \text{ for } t\in (t_0,1].
\end{equation}
By the construction of $\gamma$, $t_0$ exists.
It is sufficient to show that $t_0 = 0$.

Suppose that $t_0>0$.
Since $\rk_{\gamma(t_0)}f = n$, $f$ is locally biholomorphic at $\gamma(t_0)$.
Therefore the lifting $\psi$ can be uniquely extended
to a neighborhood $U$ of $f(\gamma(t_0))$. Set $N_1:=N\cup U$ and further
\begin{equation}
N_2 := \{ y\in N_1 : d(\psi(y),\gamma[0,1]) < d(\psi(y),f^{-1}(y_0)) \},
\end{equation}
where $d$ is the euclidean distance with respect to the ambient coordinates of $\C^n$.
Then $N_2$ and $\psi$ satisfy (3.). By passing to a smaller neighborhood
of $f(\gamma[t_0,1])$ we also obtain (1.).
Hence $t_0$ is not minimal and we have a contradiction.
\end{proof}

The following lemma is an elementary fact from linear algebra.

\begin{Lemma}\label{linear}
Let $A_1\colon V\to W_1$, $A_2\colon V\to W_2$ be linear maps
between vector spaces such that
\begin{enumerate}
\item $A_1$ is surjective;
\item $A_2|{\rm Ker} A_1$ is surjective.
\end{enumerate}
Then the direct sum $A_1\oplus A_2\colon V\to W_1\oplus W_2$
is also surjective.
\end{Lemma}

In the following lemma
let $Y_1,Y_2,Y_3$ be complex manifolds, $X\subset Z:=Y_1\times Y_2\times Y_3$
a submanifold and $a=(a_1,a_2,a_3)\in X$.
Denote by $\pi_k$ the natural projection from $X$ to
$Y_k$, $k=1,2,3$. 
Further we use the notation
\begin{equation}
X^1_a := \{z\in X : z_1=a_1\}, \quad X^2_a := \{z\in X : z_2=a_2\}
\end{equation}
and similarly $N^1_a, N^2_a$ for $N\subset Y_1\times Y_2$.

\begin{Lemma}\label{N}
Suppose that $X^1_a$ is a submanifold and 
for all $z$ from a dense subset of $X^1_a$ 
$\rk_z\pi_1 = \dim Y_1$ and $\rk_z(\pi_2|X^1_a) = \dim Y_2$.
Then there exist an open subset $N\subset Y_1\times Y_2$ 
and a holomorphic lifting
$\psi\colon N\to X$ (i.e. $(\pi_1\times\pi_2)\circ\psi=\id_N$) 
such that the following is satisfied:
\begin{enumerate}
\item $N^1_a$ is connected;
\item $(a_1,a_2)$ is in the closure of $N^1_a$;
\item $\psi$ extends continuously to $(a_1,a_2)$ 
with $\psi(a_1,a_2)=a$.
\end{enumerate} 
\end{Lemma}

\begin{proof}
Fix some local coordinates in 
$Z$ near $a$ (independently of the product structure) such that 
$X$ and $X^1_a$ equal to the unit balls in corresponding linear subspaces.
By the assumptions, there exists $z\in X^1_a$, $z_2\ne a_2$, such that
\begin{equation}
\rk_z\pi_1 = \dim Y_1, \quad \rk_z(\pi_2|X^1_a) = \dim Y_2.
\end{equation}

By Lemma~\ref{linear},
$\rk_z(\pi_1\times\pi_2)=\dim Y_1 + \dim Y_2$ and therefore the fibers 
\begin{equation}\label{100'}
(\pi_1\times\pi_2)^{-1}(w_1,w_2)
\end{equation}
through $w$ are smooth and of constant dimension for $w\in X$ close to $z$.

We claim that $z$ with this property can be chosen such that $a$ does not lie
in the tangent subspace at $z$ to the fiber $(\pi_1\times\pi_2)^{-1}(z_1,z_2)$ 
(the tangent subspace is understood as 
an affine subspace of the ambient coordinate space).
Otherwise $a$ would lie in every fiber which contradicts the assumptions.

Therefore there exists a linear subspace
$L\subset Y_1\times Y_2\times Y_3$ through $a$ and $z$ 
which is transverse and of complementary dimension to the fiber
$(\pi_1\times\pi_2)^{-1}(z_1,z_2)$. 
Then $\rk_z(\pi_1\times\pi_2)=\dim(Y_1\times Y_2)$.

In addition, we can choose $L$ transverse to $X$ and $X^1_a$ at $a$.
Since $X^1_a$ is a ball, $L\cap X^1_a$ is connected.
Hence $\rk_w(\pi_1\times\pi_2)=\dim(Y_1\times Y_2)$
for all $w$ from a dense subset of $L\cap X^1_a$.

Then $f:=(\pi_1\times\pi_2)|(L\cap X)$ 
together with $Y':=\{a_1\}\times Y_2\subset Y:=Y_1\times Y_2$
satisfies the assumptions of Lemma~\ref{eq-dim}.
The lemma follows now directly from Lemma~\ref{eq-dim}.
\end{proof}

\subsection{Applications of liftings.}

As the next step we apply the above lemmata to our situation.
Again we assume that $s:=d+1$ is odd.
Set $b:=(x,\bar x)\in\C^n\times\ov{\C^n}$.

\begin{Lemma}\label{philemma}
Under the assumptions of {\rm Theorem~\ref{main}}
for all integers $s\ge 1$, $l\ge 0$, $m:=sr+l$,
there exist an open 
subset $N\subset\C^n\times\ov{\C^n}$,
an open neighborhood $U(j_x^mf_0)\subset J^m(\C^n,\C^{n'})$
and for $(x,\bar w)\in N$, holomorphic maps
\begin{equation}
        \phi^l \colon U_x(j_x^mf_0)\times N^1_b \to J^l(\ov{\C^n},\ov{\C^{n'}}), \quad
        \Phi^l_{\bar w} \colon \ov{U_w(j_x^mf_0)}\times N^2_{(x,\bar w)} \to {\C^{n'}}
\end{equation}
such that the following is satisfied:
\begin{enumerate}
\item $N^1_b$ is connected;
\item $b$ is in the closure of $N^1_b$;
\item for all admissible germs $f$ at $x$ with
$j^m_x f\in U(j^m_x f_0)$,
\begin{equation}\label{f}
        j^l_{\bar w}\bar f = \phi^l( j^m_x f, \bar w ), \quad
        j_z^l f = \Phi^l_{\bar w}(j^m_{\bar w} \bar f, z),
\end{equation}
for $(z,\bar w)\in N$ sufficiently close to $(x,\bar x)$.
\end{enumerate}
Moreover, the map $\Phi^l_{\bar w}$
depends holomorphically on $\bar w$.
\end{Lemma}

\begin{proof}
Put in Lemma~\ref{N} 
\begin{equation}
X:=M^s\subset\C^n\times\C^N\times\ov{\C^n}, \quad
Y_1:=\C^n, \quad Y_2:=\ov{\C^n}, \quad Y_3:=\C^N
\end{equation}
and fix the point $a:=(x,\ldots,\bar x)\in X$. 
By Lemma~\ref{M-smooth}, $X$ and $X^1_a$ are manifolds and 
the first condition $\rk_z\pi_1 = \dim Y_1$ in Lemma~\ref{N}
is satisfied. The second condition $\rk_z(\pi_2|X^1_x) = \dim Y_2$
is satisfied by the choice of $s$ which we made by Theorem~\ref{opensegre}.
Let 
\begin{equation}
N\subset\C^n\times\ov{\C^n}, \quad \psi\colon N\to M^s
\end{equation} 
be given by Lemma~\ref{N}. 
Statements~(1.) and (2.) in this lemma follow directly from
(1.) and (2.) in Lemma~\ref{N} respectively.
It remains to satisfy (3.).

Let
\begin{equation}
F^l_\nu \colon U_z(j_x^mf_0) \to J^l_{\bar w}(\ov{\C^n},\ov{\C^{n'}})
\end{equation}
be as in Lemma~\ref{M-lemma}. Passing if necessary to
a smaller subset $N$ define
\begin{equation}
\phi^l(j,\bar w):=F^l_{\psi(x,\bar w)}(j)
\end{equation}
for all $j\in U_x(j_x^mf_0)$, $(x,\bar w)\in N^1_x$.
By (3.) in Lemma~\ref{N}, $\phi^l$ satisfies the 
required properties. Similarly define
\begin{equation}
\Phi^l_{\bar w}(\bar j, z):=\bar F^l_{\bar\psi(\bar w,z)}(\bar j).
\end{equation}
\end{proof}

\subsection{The end of the proof.}

The last step is based on combining the equations in (\ref{f}).
Let $f_0,f$ be as in Theorem~\ref{main} and set
\begin{equation}
m:=sr = (d+1)r, \quad k:=2m.
\end{equation}
We further use the notation of Lemma~\ref{philemma}.
The set $N^1_b$ can be seen in the canonical way
 as an open subset of $\C^n$. In general it may not contain $\bar x$.
Combining the equations in (\ref{f}) we obtain
\begin{equation}\label{12}
f(z) = \Phi^m_{\bar w}(\phi^m(j^k_x f,\bar w), z),
\end{equation}
where $(z,\bar w)\in N$ is sufficiently close to $(x,\bar x)$.

Recall that we denoted by $S$ the set of all admissible germs at $x$
(see Definition~\ref{admissible}). For the proof of Theorem~\ref{main}
we have to show conditions~(1.),(2.),(3.) in Definition~\ref{andep}.
Suppose that $j^k_xf=j^k_xf_0$. Then (1.) is followed by 
(\ref{12}).

Denote by $\Phi$ the map in the right-hand side of (\ref{12}):
\begin{equation}
\Phi\colon U_x(j^k_xf_0) \times N \to \C^{n'},
\quad
\Phi(j,z,\bar w):= \Phi^m_{\bar w}(\phi^m(j,\bar w),z).
\end{equation}
Then (\ref{12}) can be rewritten as
\begin{equation}\label{13'}
f(z) = \Phi(j^k_x f,z,\bar w).
\end{equation}
The left-hand side in (\ref{13'}) is defined for $z$ close to $x$.
The right-hand side is defined for $j^k_x f \in U(j^k_x f_0)$
and $(z,\bar w)\in N$.
To show (2.) and (3.) in Definition~\ref{andep} we fix some
$\bar w_0\in N^1_b$.

\begin{Lemma}\label{globalization}
For every admissible $f$ with $j^k_x f \in U(j^k_x f_0)$,
there exists a neighborhood $U(x)\subset\C^n$ such that
{\rm (\ref{13'})} is valid for all $z\in U(x)$ and for $w=w_0$.
\end{Lemma}

\begin{proof}
Since (\ref{13'}) is valid for all $(z,\bar w)\in N$ close to $(x,\bar x)$,
we can choose it in the form $(x,\bar w_1)$ because of (2.)
in Lemma~\ref{philemma}. Since $N^1_b$ is connected, there exists a real curve 
$\gamma\subset N^1_b$ which connects $\bar w_1$ with $\bar w_0$.
Since $\gamma$ is compact, there exists a connected product
neighborhood 
\begin{equation}
\{x\}\times\gamma \quad\subset\quad U(x)\times U(\gamma) \quad\subset\quad N.
\end{equation}
Without loss of generality, $f(z)$ is defined for all $z\in U(x)$.

Hence both sides of (\ref{13'}) are defined for $z\in U(x)$,
$\bar w\in U(\gamma)$. By the choice of $\bar w_1$, they are
equal in a neighborhood of $(x,\bar w_1)$. The conclusion of
the lemma follows from the identity principle.
\end{proof}

By Lemma~\ref{globalization}, 
for $f$ fixed and $z$ close to $x$
\begin{equation}\label{14'}
f(z) = \Phi(j^k_x f,z,\bar w_0), \quad z\in U(x).
\end{equation}
The right-hand side is defined for 
\begin{equation}
j^k_x f \in U(j^k_x f_0), \quad z\in N^2_{(x,\bar w_0)}
\end{equation}
independently of $f$.
Hence all such admissible germs $f$ extend 
to the open set $N^2_{(x,\bar w_0)}$
which proves (2.).

Finally define
\begin{equation}
F(j,z):= \Phi(j,z,\bar w_0), \quad 
F\colon U_x(j^k_xf_0)\times N^2_{(x,\bar w_0)} \to \C^{n'}.
\end{equation}
By the choice of $w_0$ one has $x\in N^2_{(x,\bar w_0)}$.
Then (3.) is implied by (\ref{14'}).
This finishes the proof of Theorem~\ref{main}.

{\bf Acknowledgments.}
I am very grateful to K.~Diederich for calling my attention to the jet method
and to M.~S.~Baouendi, P.~Ebenfelt and L.~P.~Rothschild for 
their valuable corrections and remarks.


\end{document}